\documentclass[11pt,twoside,reqno,a4paper]{amsart}

\usepackage{amsfonts,amssymb,amsmath,latexsym,epsfig,verbatim,psfrag}
\newtheorem{theorem}{Theorem}
\newtheorem{lemma}{Lemma}

\newtheorem{proposition}{Proposition}

\title[Reflected Brownian motion in a wedge]{Reflected Brownian motion in a wedge: sum-of-exponential stationary densities}
\author{A.~B.~Dieker}
\address{University College Cork, Probability Group, 17 South Bank, Crosses Green, Cork, Ireland}
\address{IBM Watson Research Center, Yorktown Heights, NY 10598}
\email{ton.dieker@isye.gatech.edu}
\author{J.~Moriarty}
\address{University College Cork, Probability Group, 17 South Bank, Crosses Green, Cork, Ireland}
\address{School of Mathematics, University of Manchester, Oxford Road, Manchester M13 9PL, UK}
\email{j.moriarty@manchester.ac.uk}

\def\proof#1{{\bf #1.}}
\def\endproof{\hfill$\Box$}
\def\sgn{\mathrm{sgn}}

\def\l{\lambda}

\def\a{\alpha}
\def\b{\beta}
\def\g{\gamma}

\def\t{\theta}
\def\d{\delta}

\def\Rot{\mathrm {Rot}}
\def\Ref{\mathrm {Ref}}
\def\tRot{\widetilde \Rot}
\def\tRef{\widetilde \Ref}

\def\R{{\mathbb R}}

\def\sgn{\mbox{sgn }}

\begin{document}

\begin{abstract}
We give necessary and sufficient conditions for the stationary density of semimartingale reflected Brownian
motion in a wedge to be written as a finite sum of terms of exponential product form. Relying on geometric ideas reminiscent of the reflection principle, we give an explicit formula for the
density in such cases. 
\end{abstract}

\maketitle

\section{Introduction}
It is well-known that Brownian motion on the positive half-line with negative drift and
reflection at zero has a stationary density which is exponential. One derivation of this fact uses time reversal to relate the distribution of the corresponding transitory process at time $t$ in the driftless case to the distribution of the maximum of a standard Brownian motion over $[0,t]$, which can be found using the reflection principle. A Cameron-Martin-Girsanov change of measure then introduces the drift, and letting $t \to \infty$ gives the required stationary distribution. In the $d$-dimensional setting, analogous arguments show that in certain situations, semimartingale reflected Brownian motion (SRBM) in a polyhedral cone has a stationary density which can be written as a finite sum of terms of exponential product form---that is, terms of the form $x
\mapsto ae^{-\langle \lambda,x\rangle}$ for some $a \in \R, \lambda \in \mathbb{R}^d$. We call such
a density a {\em sum of exponentials}.
The aim of this paper is to give necessary and sufficient conditions under which the stationary density of SRBM in a two-dimensional wedge can be written as a sum of exponentials.

In order to motivate our study, we describe the above argument in
a multidimensional setting in Section~\ref{sec:reflgroup} below.
The invariant measure for SRBM with special pushing directions at the boundary of a polyhedral cone then relates to exit
probabilities for the corresponding free Brownian motion. For some cones, in analogy with the one-dimensional case,
 these exit probabilities can be obtained explicitly using the reflection principle.
 If the reflected process is positive recurrent, then upon differentiating we obtain a stationary density
which can be written as a sum of exponentials.
Since the latter can thus be viewed as a natural multidimensional extension of the
exponential stationary density of one-dimensional SRBM, this raises the question under what
circumstances the stationary distribution of a given SRBM in a polyhedral cone has a sum-of-exponential density.
This paper shows that, for SRBM in a two-dimensional polyhedral cone,
sum-of-exponential stationary densities arise well beyond those cases to which the above arguments can be applied.

The stationary distribution of multidimensional SRBM in a polyhedral cone is only known in relatively few
cases. It has a density consisting of a single exponential term under a skew-symmetry condition on the
pushing directions on the faces \cite{harrisonwilliams:exponential1987,williams:skew1987}. We also mention an apparently
isolated example due to Harrison~\cite{harrison:tandem1978} and work of
Foschini~\cite{foschini:diffusion1982}; we discuss the latter in more detail below. In the absence
of further explicit results, numerical techniques have been developed
\cite{daiharrison:numerical1992} and logarithmic tail asymptotics have been investigated
\cite{avramdaihasenbein:explicit2001,dupuisramanan:timereversed2002,harrisonhasenbein:rbm2008}.

\begin{figure}\centering
\psfrag{x}[Br][Br]{\Large $\xi$}
\psfrag{v1}[Br][Br]{\Large $v_1$}
\psfrag{v2}[Br][Br]{\Large $v_2$}
\psfrag{n1}[Br][Br]{\Large $n_1$}
\psfrag{n2}[Br][Br]{\Large $n_2$}
\psfrag{F1}[Br][Br]{\Large $F_1$}
\psfrag{F2}[Br][Br]{\Large $F_2$}
\psfrag{tm}[Br][Br]{\Large $\theta_\mu$}
\psfrag{e}[Br][Br]{\Large $\epsilon$}
\psfrag{d}[Br][Br]{\Large $\delta$}
\psfrag{drift}[Br][Br]{\Large $-\mu$}
\resizebox{80mm}{!}{\includegraphics*{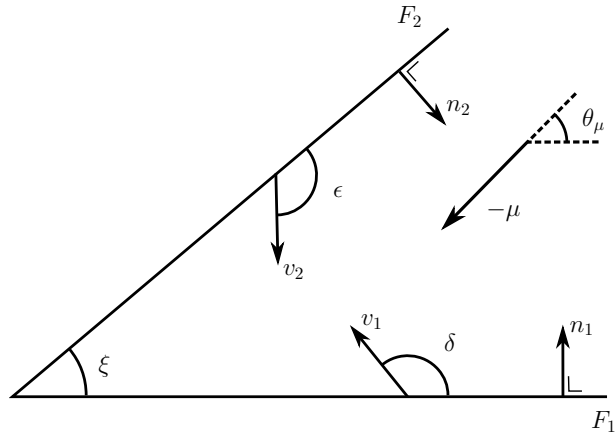}}
\caption{The wedge $S$.}
\label{fig:wedge}
\end{figure}

Before describing our results in more detail, we introduce some notation which is summarised in Figure~\ref{fig:wedge}.
Define the wedge as
\[
S=\{x\in \R^2: 0\le \arg(x)\le \xi\}.
\]
Throughout this paper we write $w_\theta=(\cos \theta,\sin \theta)$ for $\theta\in\R$. Given
some $0<\delta,\epsilon,\xi<\pi$, we set $v_1=\|v_1\| w_\delta$, $v_2=\|v_2\|
w_{\xi-\epsilon}$. Let $F_1$ and $F_2$ be the two faces of the wedge.
As is customary, we normalise $v_1$ and $v_2$ such that $\langle
v_1,n_1\rangle=1$ and $\langle v_2,n_2\rangle=1$, where $n_1$ and $n_2$ are the unit normal vectors
on $F_1$ and $F_2$, respectively. Let some vector $\mu \in \R^{2}$ also be given, and write $\theta_\mu=\arg(\mu)\in(-\pi,\pi]$.
A key role in this paper is played by $\alpha$, which is introduced by Varadhan and Williams~\cite{varadhanwilliams:wedge1985} as
\[
\alpha=\frac{\delta+\epsilon-\pi}{\xi}.
\]

Our results concern the case $\alpha<1$, and it is well-known
(e.g., \cite{varadhanwilliams:wedge1985, williams:semimartingale1985}) that under this condition there exists a
continuous semimartingale Markov process with properies that can intuitively be summarised as
follows (we do not give the mathematically rigorous definition here, as it can be found in, e.g., \cite{avramdaihasenbein:explicit2001,dupuiswilliams:lyapunov1994,williams:semimartingale1985}). The process behaves like a standard Brownian motion with drift $-\mu$ in the interior
of the wedge $S$, at the boundary of the
wedge it is pushed in some specified direction, and the time it spends at the vertex of the wedge has Lebesgue measure zero.
The pushing directions are constant along each of the faces, and are given by $v_i$ for the $i$-th face. We call this process SRBM in a wedge and seek its stationary distribution, which is absolutely continuous
with respect to Lebesgue measure \cite{dai:thesis1990,harrisonwilliams:open1987}.

We prove that the stationary density of this process is a sum of exponentials if and only if
$\alpha=-\ell$ for some nonnegative integer $\ell$.
Moreover, it follows from our results that the number of exponential terms is $2\ell+1$ when
$\alpha=-\ell$. Note that this condition reduces to $\delta=\pi-\epsilon$ if $\ell=0$, i.e.,
the skew-symmetry condition for product forms~\cite{harrisonwilliams:exponential1987}. In fact, we show that
the stationary density can be written as an $(\ell+1)\times (\ell+1)$ determinant with a special
structure.
A corollary is that the density behaves near the origin as in the driftless case studied earlier by
Williams~\cite{williams:recurrence1985}.

We next discuss some work related to our additive generalisations of the product form
exponential stationary density. To our knowledge, Foschini's study of symmetric SRBM in the
wedge~\cite{foschini:diffusion1982} was the first to consider sum-of-exponential stationary
densities. Our paper is a continuation and extension of Foschini's work, in the sense that we
couple key ideas from his algorithm with geometric arguments suggested by the reflection
principle, and consequently obtain an explicit formula for the stationary density.
A discrete-state space version of Foschini's algorithm, the {\em
compensation method}, has been successfully applied to a variety of queueing
problems \cite{adan:asymmetric1991,adan:compensation1993}.
Another body of work loosely related to the present paper studies (driftless) two-dimensional SRBM with $\alpha=-1$
\cite{dubedat:reflected2004,kager:reflected2007}, motivated by a connection with Schramm-Loewner evolutions.

This paper is organised as follows. Section~\ref{sec:reflgroup} gives some background and
discusses the special sum-of-exponential stationary densities mentioned above.
Section~\ref{sec:result} contains our main result and gives a geometric construction of the
stationary density. In Sections~\ref{sec:proofsprop}--\ref{sec:proof}, we prove the main result.

\section{Survival probabilities, reflection groups, and Weyl chambers}
\label{sec:reflgroup} This section describes a special class of sum-of-exponential stationary
densities which can be obtained directly using time reversal. These densities motivate our main result, since they typify the general structure of a sum-of-exponential stationary density. 

Consider an SRBM as in Figure~\ref{fig:wedge}, with drift $-\mu$ and
$\delta=\epsilon=\xi$. Suppose that $\mu\in S^o$ (the interior of $S$) to ensure positive
recurrence, and write $\Pi$ for the stationary measure of the SRBM. Time reversal gives
\begin{equation}
\label{eq:duality} \Pi(\{y\in S: \langle y, n_1\rangle \le \langle x, n_1\rangle, \langle y,
n_2\rangle \le \langle x, n_2\rangle\})= P_{-x}(T=\infty),
\end{equation}
where $P_y$ is the law of the corresponding free Brownian motion $B$ starting in $y$ and $T$
is the first exit time from $-S$, i.e., $T=\inf\{t\ge 0: B(t)\not \in -S\}$. To see this, we
use the wedge $-S$ to define a partial order $<$ on $\R^2$ (see for example
\cite{legall:stable1987}) and then if $\tilde{B_s}=B_{t-s}-B_t$,
\begin{eqnarray}\label{eq:rel}
\sup\{B_s:0\leq s \leq t\}-B_t&<&x \,\, \text{if and only if} \\
\sup\{-x+\tilde{B}_s:0\leq s \leq t\}&<&0, \label{eq:rel2}
\end{eqnarray}
where $\sup$ denotes the supremum with respect to the partial order. On the left hand side of
\eqref{eq:rel2} is the supremum over $[0,t]$ of the free Brownian motion starting in $-x$,
while on the left hand side of \eqref{eq:rel} we have an SRBM in $S$ as defined above (starting
in 0; for details see \cite{legall:stable1987}). Applying Wiener measure to
\eqref{eq:rel}--\eqref{eq:rel2} and letting $t \to \infty$ gives \eqref{eq:duality}.

Biane {\em et al.}~\cite{biane:littlemann2005} have recently shown
that if the wedge angle is of the form $\xi=\pi/m$ for some integer
$m\ge 2$ then
\begin{equation}
\label{eq:biane} P_{-x}(T=\infty)=\sum_{w\in G} \sgn(w) e^{-\langle \mu,(I-w)x\rangle}
\end{equation}
for any $x\in S^o$ and any $\mu\in S^o$, where the sum is taken over
the finite reflection group $G$ with associated signature function
$\text{sgn}(\cdot)$, as detailed below. We remark that this formula
gives the probability that a standard Brownian motion with drift
never exits a so-called Weyl chamber---in particular, it is not
restricted to a two-dimensional setting.

We next describe the index of summation in (\ref{eq:biane}), i.e.,
the elements of the reflection group $G$. Throughout, we represent
any element of $\R^2$ as a column vector. Write $R_\theta$ for the
reflection matrix across the line with argument $\theta$, i.e.,
\[
R_\theta=\left(\begin{array}{cc} \cos 2\theta & \sin 2\theta \\
\sin 2\theta & -\cos 2\theta \end{array}\right)
\]
and $\rho_\theta$ for the rotation matrix over an angle $\theta$, i.e.,
\[
\rho_\theta=\left(\begin{array}{cc} \cos \theta & -\sin \theta \\
\sin \theta & \cos \theta \end{array}\right).
\]
The group $G$ consists of the $2m$ matrices $I, R_\xi, \rho_{2\xi}, \rho_{2\xi}
R_\xi,\ldots,\rho_{2(m-1)\xi}, \rho_{2(m-1)\xi} R_\xi$. The signature of the matrices
$\rho_{2k\xi}R_\xi$ is $-1$, while the signature of the others is $+1$.

On combining (\ref{eq:duality}) and (\ref{eq:biane}), we find that the stationary {\em density} of SRBM with $\xi=\delta=\epsilon=\pi/m$ is proportional to
\begin{equation}
\label{eq:densityMX}
\sum_{w\in G} \sgn(w) \langle \mu, (I-w)v_2\rangle \langle \mu, (I-w)v_1\rangle e^{-\langle \mu,(I-w)x\rangle}.
\end{equation}
This is a sum of $2m-3$ nonzero terms, since the prefactor
vanishes when $w=I$, $w=R_\xi$, and $w=\rho_{2(m-1)\xi} R_\xi$.
It is the aim of this paper to put the explicit
formula (\ref{eq:densityMX}) into the context of more general SRBMs.
Specifically, we show that every sum-of-exponential stationary density has a representation
reminiscent of (\ref{eq:densityMX}).

\section{Main results}
\label{sec:result}
This section presents our main result on sum-of-exponential stationary densities, and explains its geometric interpretation. 
Our result shows that sum-of-exponential stationary densities can be written as determinants with a special structure.
The proofs are deferred to Sections~\ref{sec:proofsprop}--\ref{sec:proof}.

The matrices
\begin{equation}
\Rot_k=\rho_{2\delta+2k\xi} \quad \mathrm {and}\quad \Ref_k=\rho_{2\delta+2(k-1)\xi}R_\xi \label{eq:rotref}
\end{equation}
for $k\ge 0$ play an important role in our main theorem. Write
\[
\Theta_\ell = \{\theta \in (\xi-\epsilon,\delta): \sin(\theta-2\delta-k\xi)\neq 0 \text{ for } k=0,\ldots,2\ell\}.
\]
We know from Hobson and Rogers~\cite{hobsonrogers:rec1993} and from Dupuis and Williams~\cite{dupuiswilliams:lyapunov1994}
that there exists a unique stationary distribution of the SRBM if $\xi-\epsilon<\theta_\mu<\delta$.
Let $p^\mu$ be its density with respect to Lebesgue measure, and write $e_1$ for the vector $(1,0)$.
For integers $j$, we also define the function $\pi_j^\mu:S\to \R$ through
\[
\pi^\mu_j(x)=\frac{\langle \mu, (I-\Rot_j) v_1\rangle e^{-\langle \mu, (I-\Rot_j) x\rangle}-
\langle \mu, (I-\Ref_j) v_1\rangle e^{-\langle \mu, (I-\Ref_j) x\rangle}}{{\langle \mu, (\Ref_j-\Rot_j) v_1\rangle}}.
\]
Observe that both $p^\mu$ and $\pi_j^\mu$ depend on $\delta,\epsilon$ and $\xi$, which is suppressed
in the notation.

\begin{theorem}
\label{thm:mainresult}
The following are equivalent:
\begin{enumerate}
\item[(i)] $\alpha=-\ell$ for some integer $\ell\ge 0$;
\item[(ii)] $\alpha=-\ell$ for some integer $\ell\ge 0$ and for all $\mu$ with $\theta_\mu\in \Theta_\ell$, the functions
$p^\mu$ and $\pi^\mu$ are equal up to a multiplicative constant,
where for $x\in S$,
\begin{equation}
\pi^\mu(x) = \left|\begin{array}{cccc}
\pi^\mu_0(x) & \pi^\mu_1(x) & \cdots & \pi^\mu_\ell(x) \\
\langle \mu , \Rot_0e_1\rangle^{\ell-1} & \langle \mu , \Rot_1e_1\rangle^{\ell-1}& \cdots &\langle \mu , \Rot_\ell e_1\rangle^{\ell-1}\\
\vdots & \vdots && \vdots\\
\langle \mu , \Rot_0e_1\rangle^2 & \langle \mu , \Rot_1e_1\rangle^2 & \cdots &\langle \mu , \Rot_\ell e_1\rangle^2\\
\langle \mu , \Rot_0e_1\rangle & \langle \mu , \Rot_1e_1\rangle& \cdots &\langle \mu , \Rot_\ell e_1\rangle\\
1 & 1 & \cdots & 1
\end{array} \right|; \label{eq:determinant}
\end{equation}
\item[(iii)] $\alpha<1$ and for some $\mu$ with $\xi-\epsilon<\theta_\mu<\delta$, there exist $K<\infty$, coefficients $a_1,\ldots,a_K$,
and vectors $d_1,\ldots,d_K$ such that $p^\mu$ admits the representation
\[
p^\mu(x)= \sum_{i=1}^K a_i e^{-\langle d_i,x\rangle}.
\]
\end{enumerate}
\end{theorem}

In part (ii) of this theorem,
the condition $\theta_\mu\in \Theta_\ell$ prevents a certain degeneracy which is fundamental to the problem.
Indeed, it guarantees the linear independence of the exponential functions in $\pi^\mu$. We remark that the determinant in (ii) has a different form than the ones recently
studied in connection with transition probabilities for certain Markov chains~\cite{diekerwarren:determinant2008}.

A straightforward calculation using Vandermonde matrices shows that
the function $\pi^\mu$ defined by \eqref{eq:determinant} may be written as
\begin{equation}
\label{eq:anticlockwisepi}
\pi^\mu(x)=
\sum_{k=0}^\ell c_k \left[\langle \mu, (I-\Rot_k) v_1\rangle e^{-\langle \mu, (I-\Rot_k) x\rangle}-
\langle \mu, (I-\Ref_k) v_1\rangle e^{-\langle \mu, (I-\Ref_k) x\rangle}\right],
\end{equation}
where
\begin{equation}
\label{eq:ckvandermonde}
c_k=(-1)^k\frac{\prod_{0\le i<j\le \ell; \,i,j\neq k} \langle \mu , (\Rot_i-\Rot_j) e_1\rangle}{{\langle \mu, (\Ref_k-\Rot_k) v_1\rangle}}.
\end{equation}
This representation for $\pi^\mu$ is used later in the proof of Theorem~\ref{thm:mainresult}; note that
each exponential term is characterised by a rotation or a reflection matrix, $\Rot_0,\Ref_1,\Rot_1,\ldots,\Ref_\ell,\Rot_\ell$.

\begin{figure}\centering
\psfrag{xi}[Br][Br]{\large $\xi$}
\psfrag{x}[Br][Br]{\large $x$}
\psfrag{2d}[Br][Br]{\large $2\delta$}
\psfrag{2e}[Br][Br]{\large $2\epsilon$}
\psfrag{Refx}[Br][Br]{\small $R_\xi x$}
\psfrag{Ref0}[Br][Br]{\small $\Ref_0 x$}
\psfrag{Rot0}[Br][Br]{\small $\Rot_0 x$}
\psfrag{Ref1}[Br][Br]{\small $\Ref_1 x$}
\psfrag{Rot1}[Br][Br]{\small $\Rot_1 x$}
\psfrag{Ref2}[Br][Br]{\small $\Ref_2 x$}
\psfrag{Rot2}[Br][Br]{\small $\Rot_2 x$}
\resizebox{74mm}{!}{\includegraphics*{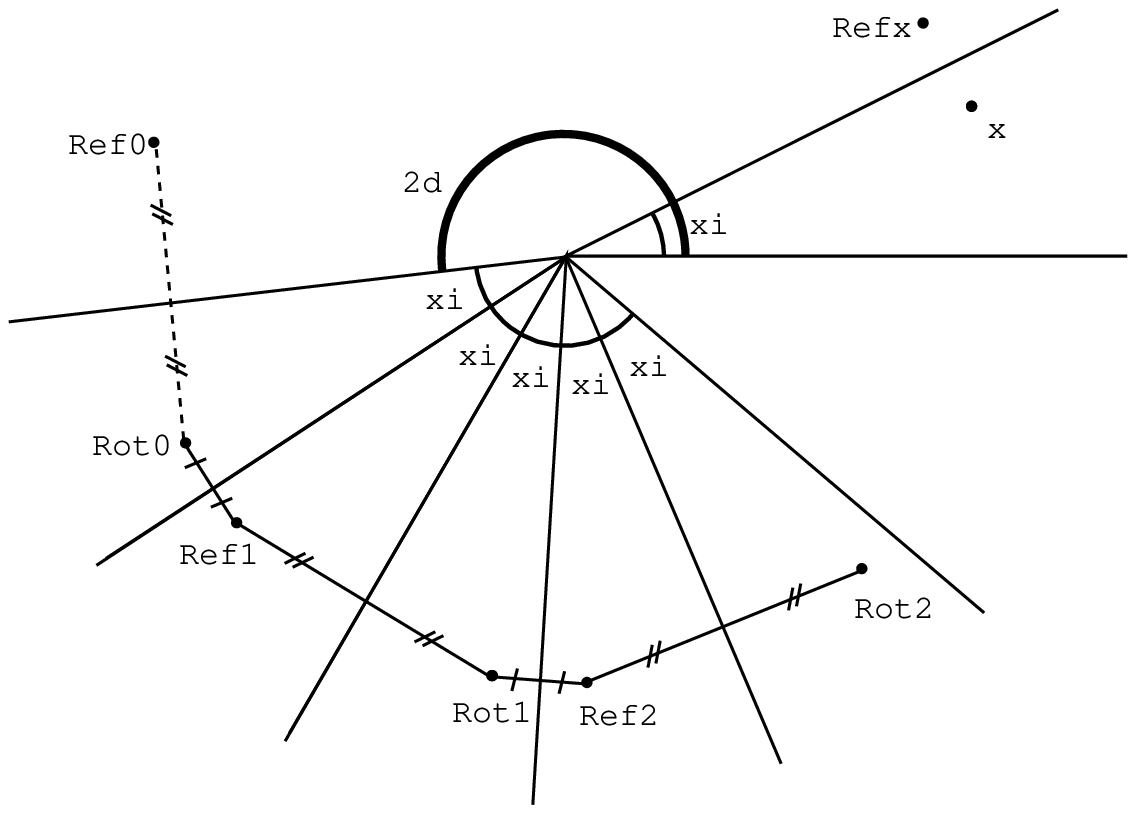}} \hspace{1mm}
\psfrag{Rot0}[Br][Br]{\small $\tRot_2 x$}
\psfrag{Ref1}[Br][Br]{\small $\tRef_2 x$}
\psfrag{Rot1}[Br][Br]{\small $\tRot_1 x$}
\psfrag{Ref2}[Br][Br]{\small $\tRef_1 x$}
\psfrag{Rot2}[Br][Br]{\small $\tRot_0 x$}
\psfrag{Ref3}[Br][Br]{\small $\tRef_0 x$}
\resizebox{74mm}{!}{\includegraphics*{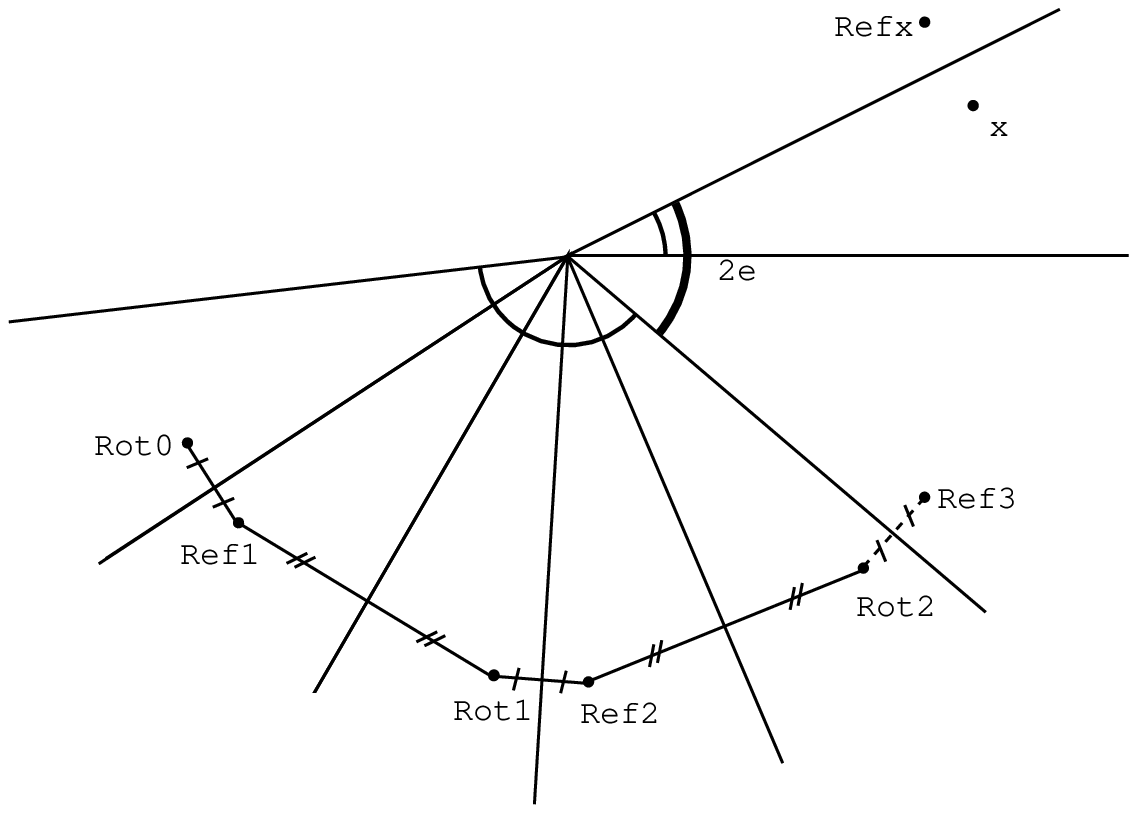}} \\
\caption{An anticlockwise construction of $\pi^\mu$ (left) and
a clockwise construction of $\tilde \pi^\mu$ (right).}
\label{fig:geometric}
\end{figure}

Although there is not necessarily an underlying reflection group in
the general setting of Theorem~\ref{thm:mainresult}, the rotation and reflection matrices in Theorem~\ref{thm:mainresult} suggest a connection with the reflection-group
framework. We illustrate this for $\alpha=-2$ (hence $\ell=2$) in the
leftmost diagram of Figure~\ref{fig:geometric}, where
we depict the five points $\Rot_0 x,\Ref_1 x,\Rot_1 x,\Ref_2 x,\Rot_2 x$
for an arbitrarily chosen $x\in S$. By construction, $\Rot_0 x$ lies
in the wedge $\rho_{2\delta} S$, which we call the `initial wedge'.
The other points $\Ref_0 x, \Rot_1 x,\ldots$ are constructed by
successive reflections, reminiscent of the orbit of $\Rot_0 x$ under
the action of a reflection group. In particular
\begin{itemize}
\item Each point lies in one of the wedges constructed by rotating the initial wedge anticlockwise over multiples of the wedge angle
$\xi$.
\item Two points lying in adjacent wedges---that is, wedges which share a common boundary line---are reflections of each
other in that line.
\end{itemize}
Note that $\Ref_0x$, although not contributing to $\pi^\mu$ since $(I-\Ref_0) v_{1}=0$, is also
obtained by reflection from $\Rot_0 x$.
We indicate this by a dashed line in Figure~\ref{fig:geometric}.

If $\alpha=-\ell$ then the last point lies in the wedge
$\rho_{2\delta+2\ell\xi} S$, which is the same as
$\rho_{-2\epsilon}S$. On comparing this with the initial wedge
$\rho_{2\delta} S$, it transpires that the {\em last} wedge in the
anticlockwise construction given above is the {\em first} wedge in
the following clockwise construction.
For $k\ge 0$ we introduce the matrices
\[
\tRot_k=\rho_{-2k\xi-2\epsilon} \quad \mathrm {and}\quad \tRef_k=\rho_{-2(k-1)\xi-2\epsilon}R_0.
\]
The rightmost diagram of Figure~\ref{fig:geometric} illustrates a clockwise construction for
the stationary density, starting with the wedge $\rho_{-2\epsilon}S$ and labelling the points using the matrices
$\tRot_k$ and $\tRef_k$. The problem is exactly the same as in the leftmost
diagram, so by uniqueness the corresponding sum-of-exponential densities must agree on $S$.
We use this observation
in the proof of Theorem~\ref{thm:mainresult}.

\medskip
We close this section by
stating some properties of the function $\pi^\mu$ defined in (\ref{eq:determinant}), which
play an important role in our proof of Theorem~\ref{thm:mainresult}. 
These properties are proved in Section~\ref{sec:proofsprop}.
We start with a result for the coefficients $\{c_k\}$ defined in (\ref{eq:ckvandermonde}). 

\begin{lemma}
\label{lem:vandermonde}
Let $-\alpha=\ell \in\{1,2,\ldots\}$ and $\theta_\mu\in\Theta_\ell$.
The coefficients $c_k$ defined by \eqref{eq:ckvandermonde}
satisfy, for $1\le k\le \ell$,
\begin{equation}
\label{eq:recck}
c_k\langle \mu, (I-\Ref_{k})v_1\rangle\langle \mu, (I-\Rot_{k-1})v_2\rangle =c_{k-1}\langle \mu, (I-\Ref_{k})v_2\rangle
\langle \mu, (I-\Rot_{k-1})v_1\rangle.
\end{equation}
\end{lemma}

Using this lemma,
it is readily checked that (\ref{eq:densityMX}) is recovered from Theorem~\ref{thm:mainresult} by setting $\delta=\epsilon=\xi=\pi/m$.
Indeed, we then have $\langle \mu, (I-\Rot_{k-1})v_1\rangle=\langle
\mu, (I-\Ref_{k})v_1\rangle$, so that (\ref{eq:recck}) reduces to
$$c_k\langle \mu, (I-\Rot_{k-1})v_2\rangle=c_{k-1}\langle \mu,
(I-\Ref_{k})v_2\rangle .$$ In this special case we may therefore set
\[
c_k=\langle \mu, (I-\Rot_k)v_2\rangle=\langle \mu, (I-\Ref_{k})v_2\rangle,
\]
and we obtain (\ref{eq:densityMX}).

The limiting behaviour of $\pi^\mu$ given in the following proposition should be compared with the invariant
measure found by Williams~\cite{williams:recurrence1985} for reflected Brownian motion {\em
without} drift. We abbreviate $\lim_{r\to0} f(r)/g(r)=1$ by $f(r)\sim g(r)$ as $r\to0$.

\begin{proposition}
\label{prop:corner} If $-\alpha=\ell\in\{1,2,\ldots\}$, then for
any $\theta\in [0,\xi]$ and any $\mu$ with $\theta_\mu\in \Theta_\ell$, we have
$\pi^\mu(rw_\theta)\sim C^\mu r^{\ell}\sin(\ell\theta+\delta)$ as $r\to 0$, where $C^\mu$ is
a finite nonzero constant independent of $r$ and $\theta$.
\end{proposition}

Our next result is that $\pi^\mu$ defined by
(\ref{eq:determinant}) does not change sign on $S$.
Note that this resolves Conjecture~1 in Dai and Harrison~\cite{daiharrison:numerical1992} for the special class of SRBMs studied in this paper.

\begin{proposition}
\label{prop:nosignchange}
The function $\pi^\mu$ does not change sign on $S$.
\end{proposition}

\section{Properties of $\pi^\mu$}
\label{sec:proofsprop}
In this section, we prove the properties of $\pi^\mu$ claimed in Section~\ref{sec:result}.
The proof of the main result, Theorem~\ref{thm:mainresult}, is deferred to Sections~\ref{sec:BAR} and \ref{sec:proof}.

\subsection{Proof of Lemma~\ref{lem:vandermonde}}
We first divide (\ref{eq:recck}) by $\sin(\theta_\mu-\delta-k\xi)\sin(\delta+(k-1)\xi)$, which
is nonzero as a consequence of the assumption on $\mu$ in conjunction with the identity
$\delta+\epsilon=\pi-\ell\xi$. Again using this identity, we find after some elementary
trigonometry that (\ref{eq:recck}) is equivalent to, with $\omega_k=\theta_\mu-2\delta-k\xi$,
\[
c_k\sin(k\xi)\sin(\omega_{\ell+k})=-c_{k-1}\sin((\ell+1-k)\xi) \sin(\omega_{k-1}).
\]
To show that this holds for the $c_k$ defined in
(\ref{eq:ckvandermonde}), we observe that
\[
\langle \mu, (\Ref_k-\Rot_k)
v_1\rangle\sin(\omega_{2(k-1)})=\langle \mu,
(\Ref_{k-1}-\Rot_{k-1}) v_1\rangle\sin(\omega_{2k})
\]
and that
\[
\langle \mu , (\Rot_i-\Rot_j) e_1\rangle=-2\sin((j-i)\xi)\sin(\omega_{i+j}).
\]
After some algebra we also find that
\begin{eqnarray*}
\lefteqn{
\sin(k\xi)\sin(\omega_{\ell+k})\sin(\omega_{2(k-1)})
\prod_{0\le i<j\le \ell; \,i,j\neq k} \langle \mu, (\Rot_i-\Rot_j) e_1\rangle}\\
&=&\sin((\ell-k+1)\xi) \sin(\omega_{k-1})
\sin(\omega_{2k})\prod_{0\le i<j\le \ell; \,i,j\neq k-1} \langle \mu , (\Rot_i-\Rot_j) e_1\rangle,
\end{eqnarray*}
and the claim follows.

\subsection{Proof of Proposition~\ref{prop:corner}}
For simplicity we suppose that $\|\mu\|=1$.
We first investigate the behaviour near zero of $\pi^\mu(x)$,
for which we rewrite $e^{\langle \mu,x\rangle}\pi^\mu(x)$ using the determinantal representation  \eqref{eq:determinant}.
A key ingredient is the identity $e^{r\cos(\eta)}=I_0(r)+2\sum_{n=1}^\infty \cos(n\eta) I_n(r)$, where $I_n$ is the modified Bessel function of the first kind.
Using this identity, after absorbing $e^{\langle\mu,x\rangle}$ into the first row, we rewrite the elements on this row as
\begin{eqnarray}
e^{\langle\mu,x\rangle}\pi^\mu_j(x) &=&
e^{\|x\|\langle \mu,\Rot_j w_\theta\rangle}-\langle \mu,(I-\Ref_j)v_1\rangle
\frac{e^{\|x\|\langle\mu,\Rot_jw_\theta\rangle}-
e^{\|x\|\langle\mu,\Ref_jw_\theta\rangle}}{\langle\mu,(\Rot_j-\Ref_j) v_1\rangle}\label{eq:col}\\
&=&I_0(\|x\|)+ 2\sum_{n=1}^\infty T_n(\cos(\omega_{2j}-\theta)) I_n(\|x\|)\nonumber\\
&&\mbox{}-2\frac{\langle \mu,(I-\Ref_j)v_1\rangle}{\langle\mu,(\Rot_j-\Ref_j) v_1\rangle}\sum_{n=1}^\infty
\left[T_n(\cos(\omega_{2j}-\theta))-T_n(\cos(\omega_{2j}+\theta))\right]I_n(\|x\|)\nonumber\\
&=&I_0(\|x\|)+ 2\sum_{n=1}^\infty T_n(\cos(\omega_{2j}-\theta)) I_n(\|x\|)\nonumber\\&&\mbox{}-2 \langle \mu,(I-\Ref_j)v_1\rangle \sum_{n=1}^\infty \sin(n\theta) U_{n-1}(\cos(\omega_{2j}))I_n(\|x\|),\nonumber
\end{eqnarray}
where we again set $\omega_k=\theta_\mu-2\delta-k\xi$, and $T_n$ and $U_n$ are the Chebyshev polynomials of the first and second kind, respectively.

In conjunction with some trigonometry, the above reasoning shows that
\[
e^{\langle \mu,x\rangle}\pi^\mu_j(x) = I_0(\|x\|)+ \frac2{\sin(\delta)}\sum_{n=1}^\infty h_{j,n}(\theta) I_n(\|x\|),
\]
where $h_{j,n}(\theta)$ is defined as
\[
\frac12 \sin(n\theta+\delta)U_n(\cos(\omega_{2j}))
-\langle \mu,v_1/\|v_1\|\rangle \sin(n\theta)U_{n-1}(\cos(\omega_{2j}))
+\frac 12 \sin(n\theta-\delta)U_{n-2}(\cos(\omega_{2j})).
\]
(We use the convention $U_{-1}(x)=0$.)
Therefore, $e^{\langle \mu,x\rangle}\pi^\mu(x)$ can be expanded in terms of modified Bessel functions of the first kind, and for $n\ge 1$ the coefficient in front of $I_n(\|x\|)$ is proportional to
\begin{equation}
\label{eq:chebdet}
\left|\begin{array}{ccccc}
h_{0,n}(\theta) & h_{1,n}(\theta) &h_{2,n}(\theta)& \cdots & h_{\ell,n}(\theta) \\
U_{\ell-1}(\cos(\omega_{0})) & U_{\ell-1}(\cos(\omega_{2})) & U_{\ell-1}(\cos(\omega_{4})) &  \cdots &U_{\ell-1}(\cos(\omega_{2\ell}))\\
\vdots & \vdots && \vdots\\
U_0(\cos(\omega_{0})) & U_0(\cos(\omega_{2}))& U_0(\cos(\omega_{4})) &\cdots &U_0(\cos(\omega_{2\ell}))
\end{array} \right|.
\end{equation}
To see how this follows from~\eqref{eq:determinant},
note that we may apply elementary determinantal operations to replace a row with elements $\cos(\omega_{2j})^m$ by $U_m(\cos(\omega_{2j}))$.
The term $I_0(\|x\|)$ is not present in the expansion in view of the last row in (\ref{eq:chebdet}) with ones.

The condition $\theta_\mu\in \Theta_\ell$ guarantees that none of the $\cos(\omega_{2j})$ are equal, and
we conclude that the coefficient of $I_n(\|x\|)$ vanishes for $n< \ell$
and that it is proportional to $\sin(\ell \theta+\delta)$ for $n=\ell$.
Since $I_\ell(r)\sim Cr^{\ell}$ for some constant $C\neq 0$ as $r\to0$, this yields
$\pi^\mu(rw_\theta)\sim C^\mu r^{\ell}\sin(\ell\theta+\delta)$.

\subsection{Proof of Proposition~\ref{prop:nosignchange}}
We rely on the following auxiliary lemma, whose proof is inspired by some
elementary symmetric-function theory. Alternatively, as
communicated to us by Sean Meyn, Proposition~\ref{prop:nosignchange} can be proved using a continuous-space
analogue of Theorem~1 of Foster~\cite{foster:stochastic1953}; this can be
derived using the general theory of Markov processes.

\begin{lemma}
\label{lem:signschur}
Let $\ell\in\{1,2,\ldots\}$ and $\zeta\in\R^{\ell+1}$. For any $y> 0$, the sign of the determinant
\[
\left|\begin{array}{cccc}
e^{\zeta_0 y} & e^{\zeta_1 y} & \cdots & e^{\zeta_\ell y} \\
\zeta_0^{\ell-1} & \zeta_1^{\ell-1}& \cdots &\zeta_\ell^{\ell-1}\\
\vdots & \vdots && \vdots\\
\zeta_0 & \zeta_1 & \cdots &\zeta_\ell\\
1 & 1 & \cdots & 1
\end{array} \right|
\]
equals the sign of $\prod_{0\le i<j\le \ell}[\zeta_i-\zeta_j]$.
\end{lemma}
\proof{Proof} The statement is a continuous analogue of the claim
that $s_{(n,0,\ldots,0)}(\zeta)$ is nonnegative for $\zeta\ge 0$,
where $s_\lambda$ is a symmetric polynomial known as a Schur
polynomial (or, in this special case, a complete homogeneous
symmetric polynomial).

By induction on $\ell$ one can show that the given determinant equals
\[
\prod_{0\le i<j\le \ell}[\zeta_i-\zeta_j] \int_{0= z_{-1} \le
z_{0}\le \ldots\le z_{\ell-1}\le z_\ell= y}
e^{\zeta_{0}(z_{0}-z_{-1})} \cdots
e^{\zeta_\ell(z_\ell-z_{\ell-1})} dz_{0}\cdots dz_{\ell-1},
\]
and the claim follows.
\endproof

\medskip
We now prove Proposition~\ref{prop:nosignchange}.
By the Maximum Principle (see Theorem 2.5 of
\cite{protterweinberger:maximumprin1967} for a suitable form),
neither the minimum nor the maximum of $\pi^\mu$ over $S$ is attained in
the open set $S^o$. We therefore investigate the boundary.

We first prove that $\pi^\mu(x)\to 0$ as $\|x\|\to \infty$ by showing
that, for any $x\in S$,
\begin{eqnarray*}
\langle \mu, (I-\Rot_k) x\rangle>0, &&k=0,\ldots,\ell, \\
\langle \mu, (I-\Ref_k) x\rangle>0, &&k=1,\ldots,\ell.
\end{eqnarray*}
Set $\theta=\arg x$. For the claim involving $\Rot_k$, we observe that
\[
\langle \mu, (I-\Rot_k) x\rangle=
-2\|x\|\sin(\delta+k\xi)\sin(\theta_\mu-\theta-\delta-k\xi).
\]
Since $0<\delta+k\xi\le \delta+\ell\xi=\pi-\epsilon<\pi$ and we have
$\theta < \xi$ and the stability condition
$\xi-\epsilon<\theta_\mu<\delta$, we obtain
$-\pi<\theta_\mu-\theta-\delta-k\xi<0$. The same argument works for
the claim involving $\Ref_k$, now relying on
\[
\langle \mu, (I-\Ref_k)
x\rangle=-2\|x\|\sin(\delta+k\xi-\theta)\sin(\theta_\mu-\delta-k\xi)
\]
and the assumption $k\ge 1$.

We next prove that the signs of $\pi^\mu(rw_0)$ and $\pi^\mu(rw_\xi)$ are
equal and independent of $r>0$. The equality of the signs follows
from Proposition~\ref{prop:corner} after showing that they do not
depend on $r$. From Lemma~\ref{lem:signschur} with $\zeta_j=\langle
\mu,  \Rot_j e_1\rangle$ and (\ref{eq:col}) we conclude this for
$\pi^\mu(rw_0)$. Applying the same argument `clockwise' shows that this
also holds for $\pi^\mu(rw_\xi)$.

\section{The BAR and a PDE with boundary conditions}
\label{sec:BAR}
This section prepares for the proof of Theorem~\ref{thm:mainresult} by relating the stationary
density to a partial differential equation (PDE) with boundary conditions involving the pushing directions.

\subsection{The BAR}
Our proof of Theorem~\ref{thm:mainresult} requires the Basic Adjoint Relationship (BAR) as
presented in following proposition, which is implied by Propositions~3 and 4 in
\cite{daiharrison:numerical1992}; see \cite{dai:thesis1990} for proofs.

\begin{proposition}
\label{prop:BAR} Suppose that $\alpha<1$ and 
assume the existence and uniqueness of a stationary distribution for the SRBM.

A nonzero finite measure $\nu_0$ on $S$ is proportional to the stationary distribution
if and only if
there exist finite measures $\nu_1$ on $F_1$ and $\nu_2$ on $F_2$ such that for any $f\in C^2_b(S)$
\[
\int_S \left[\frac 12 \Delta f -\langle \mu, \nabla f\rangle\right] d\nu_0 + \int_{F_1} \langle v_1,\nabla f\rangle d\nu_1+
\int_{F_2} \langle v_2,\nabla f\rangle d\nu_2=0
\quad \quad\quad\rm{(BAR)}.
\]
\end{proposition}

Let $\sigma$ and $\sigma_i$ be the Lebesgue measures on $S$ and $F_i$, respectively.
Write $v_i^*=2n_i-v_i$.

\begin{proposition}
\label{prop:partdiffeq}
Let $p\in C^2(S)$ be nonnegative and integrable over $S$.

If (BAR) is satisfied with $d\nu_0=p \, d\sigma$ and $d\nu_i=p/2 \, d\sigma_i$, then
\begin{eqnarray}
    \Delta p + 2\langle \mu, \nabla p\rangle&=&0 \quad \mbox{ on $S^o$, }\label{pde}\\
    \langle v_1^*,\nabla p\rangle+2\langle\mu, n_1\rangle p&=&0 \quad \mbox{ on $F_1^o$, }\label{bc1}\\
    \langle v_2^*,\nabla p\rangle+2\langle\mu, n_2\rangle p&=&0 \quad \mbox{ on $F_2^o$. }\label{bc2}
\end{eqnarray}

Conversely, if (\ref{pde})--(\ref{bc2}) hold and moreover $p(0)=0$, then
(BAR) is satisfied with $d\nu_0=p \, d\sigma$ and $d\nu_i=p/2 \, d\sigma_i$.
\end{proposition}
\proof{Proof}
We may repeat the arguments in the proof of Lemma~7.1 of Harrison and Williams~\cite{harrisonwilliams:exponential1987}.
The additional assumption $p(0)=0$ ensures that (7.8) in \cite{harrisonwilliams:exponential1987} automatically holds.
\endproof

\medskip
The above proposition motivates investigating sum-of-exponential solutions to (\ref{pde})--(\ref{bc2}).

\subsection{Sum-of-exponential solutions to the PDE}
In this subsection we study some properties of sum-of-exponential solutions to the PDE
(\ref{pde}) plus either boundary condition (\ref{bc1}) or (\ref{bc2}). We use the following
observation, due to Foschini~\cite[Sec.~III.A]{foschini:diffusion1982}.

\begin{lemma}
\label{lem:foschini} Let $p$ be given by $p(x)=\sum_{i=1}^k a_i e^{-\langle c_i,x\rangle}$ for
some $k<\infty$, $a_i\neq 0$, and $c_i\neq 0$ such that $c_i\neq c_j$ if $i\neq j$.

If $p$ satisfies (\ref{pde}) and (\ref{bc1}), then for each $i=1,\ldots,k$ precisely one of
the following holds:
\begin{enumerate}
\item $x\mapsto e^{-\langle c_i,x\rangle}$ satisfies (\ref{bc1}), or
\item there exists a unique $j\neq i$ such that
$x\mapsto a_i e^{-\langle c_i,x\rangle}+a_j e^{-\langle c_j,x\rangle}$ satisfies (\ref{bc1}) and
we have $\langle c_i,w_0\rangle=\langle c_j,w_0\rangle$.
\end{enumerate}
\end{lemma}

By symmetry, Lemma~\ref{lem:foschini} also holds when (\ref{bc1}) is replaced by (\ref{bc2}),
provided the condition $\langle c_i,w_0\rangle=\langle c_j,w_0\rangle$ is replaced by $\langle
c_i,w_\xi\rangle=\langle c_j,w_\xi\rangle$. The next two lemmas investigate the two scenarios
of Lemma~\ref{lem:foschini} in more detail; Lemma \ref{lem:singleexp} may be regarded as a
generalisation of Theorem 6.1 of \cite{harrisonwilliams:exponential1987} (modulo the discussion
of the BAR in Section~\ref{sec:proof}).

\begin{lemma}
\label{lem:singleexp} Let $p(x)= e^{-\langle c,x\rangle}$ for some $c\neq 0$.
\begin{enumerate}
\item
If $p$ satisfies (\ref{pde}) and (\ref{bc1}), then either $p(x)=e^{-\langle\mu,(I-\rho_{2\d})x\rangle}$ or $p(x)=e^{-\langle\mu,(I-R_{0})x\rangle}$.
\item
If $p$ satisfies (\ref{pde}) and (\ref{bc2}), then either $p(x)=e^{-\langle\mu,(I-\rho_{-2\epsilon})x\rangle}$ or $p(x)=e^{-\langle\mu,(I-R_\xi)x\rangle}$.
\end{enumerate}
\end{lemma}
\proof{Proof} We only prove the first claim, the second being the clockwise analogue. The
condition that $p$ satisfies (\ref{pde}) translates to $\|c-\mu\|=\|\mu\|$. We may therefore
write $c=\mu-\rho_{-2\gamma}\mu$ for some $\gamma\equiv \gamma(\mu)$. Next we substitute this
in (\ref{bc1}), yielding $-\langle \mu,(I-\rho_{2\gamma})v_1^*\rangle+2\langle
\mu,n_1\rangle=0$, which we may rewrite as $\langle \mu, v_1+\rho_{2\gamma}v_1^*\rangle=0$
since $v_1^*=2n_1-v_1$. Using $v_1^*=-R_\delta \rho_{2\delta} v_1$, we get $\langle \mu,
(I-\rho_{2(\gamma-\delta)}) v_1\rangle=0$. This can only be the case if $\gamma=\delta\mod
\pi$ or $\gamma=\theta_\mu\mod \pi$.
\endproof

\medskip
The next result investigates the second scenario of Lemma \ref{lem:foschini}.

\begin{lemma}
\label{lem:pair}
Let $p$ be given by $p(x)=a_1 e^{-\langle c,x\rangle}+ a_2 e^{-\langle d, x\rangle}$ for some $a_1,a_2\neq0$, $c,d\neq 0$, and $c\neq d$.
\begin{enumerate}
\item
If $\langle c,w_0\rangle=\langle d,w_0\rangle$ and $p$ satisfies (\ref{pde}) and (\ref{bc1}),
then there exists some $\gamma\equiv\gamma(\mu)\in(0,\pi)$ such that $p$ is proportional to $p_\gamma$ defined by
\[
p_\gamma(x)=
\langle \mu, (I-\rho_{2\gamma+2\delta}) v_1\rangle e^{-\langle \mu, (I-\rho_{2\gamma+2\delta}) x\rangle}-
\langle \mu, (I-\rho_{2\gamma+2\delta} R_0) v_1\rangle e^{-\langle \mu, (I-\rho_{2\gamma+2\delta} R_0)x\rangle}.
\]
\item
If $\langle c,w_\xi\rangle=\langle d,w_\xi\rangle$ and $p$ satisfies (\ref{pde}) and (\ref{bc2}),
then there exists some $\tilde \gamma\equiv\tilde \gamma(\mu)\in(0,\pi)$ such that $p$ is proportional to $\tilde p_{\tilde \gamma}$ defined by 
\[
\tilde p_{\tilde \gamma}(x)=
\langle \mu, (I-\rho_{-2\tilde \gamma-2\epsilon}) v_2\rangle e^{-\langle \mu, (I-\rho_{-2\tilde \gamma-2\epsilon}) x\rangle}-
\langle \mu, (I-\rho_{-2\tilde \gamma-2\epsilon} R_{\xi})v_2\rangle e^{-\langle \mu, (I-\rho_{-2\tilde \gamma-2\epsilon} R_{\xi})x\rangle}.
\]
\end{enumerate}
\end{lemma}
\proof{Proof} Again we only prove the first claim. By linear independence both $e^{-\langle
c,x\rangle}$ and $e^{-\langle d,x\rangle}$ must satisfy (\ref{pde}) individually. As in the
proof of Lemma~\ref{lem:singleexp}, we may therefore write $c=\mu-\rho_{-2\delta-2\gamma}\mu$
for some $\gamma\equiv \gamma(\mu)\in[0,\pi)$, so that $\langle
c,x\rangle=\langle\mu,(I-\rho_{2\gamma+2\delta}) x\rangle$. From $\langle c,w_0\rangle=\langle
d,w_0\rangle$ and $c\neq d$ we conclude that 
$\langle d, x\rangle=
\langle\mu,(I-\rho_{2\gamma+2\delta} R_0) x\rangle.$

It remains to study $a_1$ and $a_2$, for which we use (\ref{bc1}). Since $\langle
c,w_0\rangle= \langle d,w_0\rangle$
we obtain that on $F_1$,
\[
\langle v_1^*,\nabla p(x)\rangle = -[a_1\langle \mu, (I-\rho_{2\gamma+2\delta}) v_1^*\rangle +
a_2\langle \mu, (I-\rho_{2\gamma+2\delta} R_0) v_1^*\rangle]e^{-\langle c, x\rangle}.
\]
With $v_1^*=2n_1-v_1$, we conclude that (\ref{bc1}) implies
\[
a_1 \langle\mu, v_1+\rho_{2\gamma+2\delta} v_1^*\rangle+a_2 \langle\mu, v_1+\rho_{2\gamma+2\delta} R_0 v_1^*\rangle=0.
\]

The result follows after using $v_1^*=-R_0 v_1$ and noting that
$\langle \mu, (I-\rho_{2\gamma+2\delta} R_0) v_1\rangle$ cannot be zero.
\endproof

\medskip
We remark that Lemmas~\ref{lem:foschini}--\ref{lem:pair} show that some structure from the particular example of Section~\ref{sec:reflgroup} holds in general. Specifically, each exponent in a sum-of-exponential solution equals $-\langle\mu,(I-M)x\rangle$ for some reflection or rotation matrix $M$.

\section{Proof of Theorem~\ref{thm:mainresult}} \label{sec:proof}
To prove our main result, it suffices to show that (i) implies (ii) and that (iii) implies (i).
We start with the latter.

\subsection{Proof that (iii) implies (i)}
By assumption, there is some $\mu$ with $\xi-\epsilon<\theta_\mu<\delta$ such that
the stationary density has the form
\begin{equation}
\label{eq:pfinitesum}
p^\mu(x)=\sum_{i=1}^k a_i e^{-\langle c_i,x\rangle}.
\end{equation}
We assume without loss of generality that the $a_i$ are nonzero and that the $c_i$ are
distinct. We may also restrict ourselves to the case of more than one summand ($k>1$), since
the $k=1$ case has already been studied \cite{harrisonreiman:on1981,harrisonwilliams:exponential1987}.

We next argue that $p^\mu$ satisfies (\ref{pde})--(\ref{bc2}). To this end, in view of
Proposition~\ref{prop:partdiffeq}, we need to show that (BAR) is satisfied with
$d\nu_0=p^\mu \, d\sigma$ and $d\nu_i=p^\mu/2 \, d\sigma_i$. We do so using an argument due to Harrison and
Williams~\cite[p.~108]{harrisonwilliams:open1987}. From Proposition~\ref{prop:BAR} we know
that (BAR) holds for some measures $\nu_1,\nu_2$. Let $\lambda\in\R^2$ satisfy $\langle \l, x
\rangle\ge 0$ for all $x\in S$. On substituting $f(x)=e^{-\langle \lambda,x\rangle}$ in (BAR)
we find
\[
\left[\frac 12 \|\lambda\|^2+\langle \mu,\lambda\rangle\right] \int_S e^{-\langle \lambda,x\rangle}
p^\mu(x)dx-\langle v_1,\lambda\rangle \int_{F_1} e^{-\langle \lambda,x\rangle}\nu_1(dx)-\langle v_2,\lambda\rangle
\int_{F_2} e^{-\langle \lambda,x\rangle}\nu_2(dx)=0.
\]
Let us first focus on $\nu_1$. Write $\lambda=\lambda_1 w_0+\lambda_2 n_1$, and let $\lambda_2\to\infty$
after dividing the above equation by $\lambda_2$.
To evaluate the resulting limit of the first term, we use the initial value theorem
to obtain
\[
\lim_{\lambda_2\to\infty} \lambda_2 \int_S e^{-\langle \lambda,x\rangle} p^\mu(x)dx=\int_0^\infty e^{-\lambda_1 s} p^\mu(sw_0) ds.
\]
After taking the limits of all other terms as well (recalling that $\langle v_1, n_1\rangle=1$),
we conclude that $\int_{F_1} e^{-\lambda_1 x_1}\nu_1(dx)=\int_0^\infty e^{-\lambda_1 z}p^\mu((z,0))
dz/2$ for $\lambda_1\ge 0$. The uniqueness
theorem for Laplace transforms thus yields $d\nu_1=p^\mu/2 \, d\sigma_1$. A similar argument works
to show $d\nu_2=p^\mu/2 \, d\sigma_2$ by studying $\lambda=\lambda_1 w_\xi+\lambda_2 n_2$ for large
$\lambda_2$.

In the remainder of this subsection, it is our aim to further specify the structure of $p^\mu$ defined in (\ref{eq:pfinitesum}) when it solves (\ref{pde})--(\ref{bc2}). Lemmas~\ref{lem:foschini}--\ref{lem:pair} play a central role in this analysis.

\subsection*{Graph representation.} It is convenient to represent $p^\mu$ by an undirected labelled graph $G$, with $k$ vertices, as follows. Each
vertex represents a summand $a_ie^{-\langle c_i,x \rangle}$ in (\ref{eq:pfinitesum}), and we
say that two vertices are joined by a {\em BC1 edge} (respectively {\em BC2 edge}) if the sum
of the terms corresponding to these vertices satisfies \eqref{pde} and \eqref{bc1}
(respectively \eqref{pde} and \eqref{bc2}). Note that by Lemma \ref{lem:foschini}, at most one
BC1 edge and at most one BC2 edge can be incident at any given vertex. Therefore, the degree
of the vertices in $G$ cannot exceed two, and BC1 edges and BC2 edges alternate along any
path. Here and throughout, we say that a subgraph of $G$ is a {\em path} if it is connected
and acyclic. The \emph{length} of a path equals its number of vertices.
The vertex corresponding to the summand $e^{-\langle c,x\rangle}$ is labelled by the
matrix $M$ for which $\langle c,x \rangle =\langle \mu, (I-M)x\rangle$, and we refer to this vertex as an $M$ vertex. Note that $M$ is necessarily a reflection or rotation matrix. Although $M$ is uniquely defined for any of the vertices of $G$, it has several representations---for example, we also refer to an $R_\xi$ vertex as a $\rho_{2\xi} R_0$ vertex. We
refer to a path between an $M_1$ vertex and and $M_2$ vertex as an `$M_1-M_2$ path'.

\subsection*{Mating procedure} Given the label of any vertex in $G$, we can specify the labels of all other vertices
in the same connected component. For let $M_{1}$, $M_{2}$ be the labels of two arbitrary vertices joined by a BC1 edge: then by Lemma~\ref{lem:pair}, the sum of the
corresponding exponential terms is proportional to $p_\gamma$ for
some $\gamma\in(0,\pi)$. We will say that $M_{1}$ is the \emph{BC1 mate} of $M_{2}$.
By considering separately the cases when $M_{1}$ is a reflection
and a rotation, it is easy to see that $\{M_{1},M_{2}\}=\{\rho_{2\beta},\rho_{2\beta-2\xi}R_\xi\}$ for some angle
$\beta\in(-\pi,\pi]$.
Similarly, the labels of two arbitrary vertices joined by a BC2 edge are $\{\rho_{2\beta},\rho_{2\beta}R_\xi\}$ for some $\beta\in(-\pi,\pi]$.
Any path in $G$---beginning for example with a BC2 edge---therefore has labels $\rho_{2\beta},\rho_{2\beta}R_\xi,\rho_{2\beta+2\xi},\ldots$ for some $\beta$.

\subsection*{Example.} 
\begin{figure}\centering
\psfrag{BC1}[Bc][Bc]{BC1}
\psfrag{BC2}[Bc][Bc]{BC2}
\psfrag{J1}[Bc][Bc]{\large $\rho_{2\delta}$}
\psfrag{J2}[Bc][Bc]{\large $\rho_{2\delta}R_{\xi}$}
\psfrag{J3}[Bc][Bc]{\large $\rho_{2\delta+2\xi}$}
\psfrag{J4}[Bc][Bc]{\large $\rho_{2\delta+2\xi}R_{\xi}$}
\psfrag{J5}[Bc][Bc]{\large $\rho_{2\delta+4\xi}$}
\resizebox{90mm}{!}{\includegraphics*{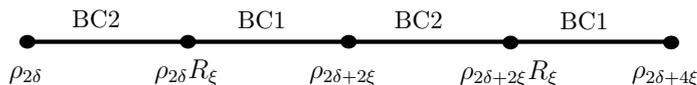}} \hspace{1mm}
\caption{An example labelled graph. By \eqref{eq:rotref}, this graph corresponds to the reflection construction in the leftmost diagram of Figure~\ref{fig:geometric}. If $\a=-2$ then this graph also corresponds to the rightmost diagram of Figure~\ref{fig:geometric}.}
\label{fig:ghost}
\end{figure}

To illustrate the graph representation and mating procedure, suppose that one summand in \eqref{eq:pfinitesum} has  exponent $-\langle\mu,(I-\rho_{2\delta})x\rangle$. In the graph representation, this summand is represented by a vertex with label $\rho_{2\d}$: suppose there exists a path of length 5 starting at this vertex, and that its first edge is a BC2 edge.
By the mating procedure, the vertex labels for this path are $\rho_{2\delta},\rho_{2\delta}R_{\xi},\rho_{2\delta+2\xi},\rho_{2\delta+2\xi}R_{\xi},\rho_{2\delta+4\xi}$, cf.~Figure~\ref{fig:ghost}. By \eqref{eq:rotref}, this path corresponds to the reflection construction in the leftmost diagram of Figure~\ref{fig:geometric}.

\begin{proposition} \label{prop:item} Let the stationary density $p^\mu$ be of the form \eqref{eq:pfinitesum}. If $G$ is the labelled graph corresponding to $p^\mu$ then:
\begin{itemize}
\item $G$ is a $\rho_{2\delta}-\rho_{-2\epsilon}$ path;
\item The number of vertices in $G$ is odd;
\item $\alpha=-\ell$ for some integer $\ell\ge 0$.
\end{itemize}
\end{proposition}
\proof{Proof} As usual we exclude the case $\a=0$. Since $p^\mu$ is a density, each exponent in \eqref{eq:pfinitesum} is nonpositive and $p^\mu(x)\to 0$ as $\|x\|\to\infty$ in $S$.

A key tool in the proof is the following {\em range restriction} for the reflection labels.
That is, $G$ cannot contain an $R_\gamma$ vertex if $w_\gamma\in S^o$ or $-w_{\gamma}\in S^o$.
Suppose {\sl a contrario} that $G$ contains an $R_\gamma$ vertex and that $L_{\gamma}\cap S^o\not =\emptyset$,
where $L_{\gamma}$ is the line $\{rw_{\gamma}:r\in\R\}$.
We must then have $\mu\in L_\gamma$ since otherwise the exponent $x\mapsto -\langle \mu, (I-R_{\gamma})x \rangle$ changes sign in $S^o$ on either side of $L_\gamma$. However, if $\mu\in L_\gamma$ then $-\langle \mu, (I-R_{\gamma})x \rangle = 0$ for all $x$, but such a constant exponent (which must be unique) contradicts $p^\mu(x)\to0$.

\emph{$G$ is acyclic.} Suppose that $G$ contains a cycle $G_{0}$, of length $2m$ say (note that cycles of odd length are impossible by the mating procedure).
Taking an arbitrary rotation label $\rho_{2\b}$ from $G_{0}$, we must have $\rho_{2\b}=\rho_{2\b + 2m\xi}$ by cyclicity and the mating procedure, so that $\xi=n\pi/m$ for some integer $n \ge1$.
The reflection labels in $G_0$ are readily seen to be $R_{\beta+k\xi}, k = 1,\ldots,m$, so
range restriction yields $L_{\beta+k\xi}\cap S^o=\emptyset$ for each $k=1,\ldots,m$. Using this in conjunction with the fact that (by uniqueness of the reflection labels) $\{L_{\beta+kn\pi/m}:k=1,\ldots,m\}=\{L_{\beta+k\pi/m}: k=1,\ldots,m\}$, we deduce that $n=1$ and that $R_{\beta+k\xi}=R_0$ for some $k$. By the proof of Lemma \ref{lem:singleexp} $R_0$ must be a vertex of degree 1, so this is a contradiction.

\emph{$G$ does not contain the label $R_0$, nor does it contain the label $R_\xi$.} We only prove that $G$ does not contain the vertex $R_0$; similar arguments can be given for $R_\xi$.
Suppose that $G$ contains the vertex $R_0$, and consider the sum-of-exponentials $p^\mu(x)$
corresponding to $G$ when $x$ lies on the boundary $F_1$.
There is a constant nonzero term since $(I-R_0)w_0=0$, so to ensure $p^\mu(rw_0)\to 0$ as $r\to\infty$ there must be another exponent that vanishes on $F_1$. First observe that $R_0$ is the only reflection label that can possibly correspond to a constant term on $F_1$ (as shown already, $\mu\not \in L_\gamma$ if $G$ contains the label $R_\gamma$ so $\mu$ cannot be orthogonal to $(I-R_\g)w_0$).
All rotation vertices in $G$ correspond to exactly the same exponents on $F_1$ as their possible BC1 mates, which are reflection vertices necessarily different from $R_0$, whence none of the vertices joined by a BC1 edge can represent the term constant on $F_1$. Therefore, in view of Lemma~\ref{lem:singleexp}, the only remaining possibility is that $G$ contain a vertex labelled $\rho_{2\delta}$ which corresponds to the constant on $F_1$. However, this contradicts our assumption that $\xi-\epsilon<\theta_\mu<\delta$.

We have now proved that $G$ is a $\rho_{2\delta}-\rho_{-2\epsilon}$ path, and by the mating procedure the number of terms is odd, say $2\ell+1$. The mating procedure also shows that we must have $\rho_{-2\epsilon}= \rho_{2\ell\xi + 2\delta}$, so that $\delta+\epsilon+\ell\xi=n \pi$ for some integer $n \ge 1$. The range restriction on the reflection labels entails that $L_{\delta+\xi},L_{\delta+2\xi},\ldots,L_{\delta+\ell\xi}\not \in S^o$. None of these lines can be equal to $L_0$, since this would imply that $R_0$ is a label in $G$. We deduce that $\delta+\ell\xi<\pi$, and therefore $n=1$.~\endproof

\medskip
To continue our example, suppose that $\a=-2$. Then $\rho_{2\delta+4\xi}=\rho_{-2\epsilon}$ and so the graph in our example also corresponds to the rightmost diagram in Figure~\ref{fig:geometric}. Then $G$ has 5 vertices, so that $p^\mu$ has $k=5$ exponential terms.

\subsection{Proof that (i) implies (ii)}
Suppose that $\alpha=-\ell$ for some $\ell\in\{1,2,\ldots\}$, and consider a $\mu$ with
$\theta_\mu\in \Theta_\ell$. We shall use the representation \eqref{eq:anticlockwisepi} and Lemma~\ref{lem:vandermonde} to argue that $\pi^\mu$
must equal the stationary density up to a multiplicative constant. 
We first argue that $\pi^\mu$ satisfies (\ref{pde})--(\ref{bc2}).

The proofs of Lemmas~\ref{lem:singleexp} and \ref{lem:pair} show that $\pi^\mu$
satisfies (\ref{pde}) and (\ref{bc1}). To see that $\pi^\mu$ also satisfies (\ref{bc2})
we note that, for any constants $d_0,\ldots,d_\ell$, the
function $\tilde \pi^\mu$ defined on $S$ by
\begin{equation}
\label{eq:tildepid}
\tilde \pi^\mu(x)=
d_0 \langle\mu,(I-\rho_{-2\epsilon}) v_2 \rangle e^{-\langle\mu,(I-\rho_{-2\epsilon}) x\rangle} + d_1 \tilde p_\xi(x) + d_2 \tilde p_{2\xi}(x)+\ldots+ d_\ell \tilde p_{\ell\xi}(x)
\end{equation}
satisfies (\ref{pde}) and (\ref{bc2}) ($\tilde p_\theta$ is defined in Lemma \ref{lem:pair}). 
Using $\epsilon+\delta+\ell\xi=\pi$, on investigating the exponents we find that $\pi^\mu$ and $\tilde\pi^\mu$ are linear combinations of the same exponential terms, so it suffices to show that the coefficients are proportional to
each other. To do so, we write
\begin{eqnarray*}
\pi^\mu(x)&=& c_\ell \langle \mu, (I-\Rot_\ell) v_1\rangle e^{-\langle \mu, (I-\Rot_\ell) x\rangle}\\
&&\mbox{}+
\sum_{k=1}^\ell \left[c_{\ell-k}
\langle \mu, (I-\Rot_{\ell-k}) v_1\rangle e^{-\langle \mu, (I-\Rot_{\ell-k}) x\rangle}\right.\\
&&\left.\hspace{12mm}\mbox{}-c_{\ell-k+1}
\langle \mu, (I-\Ref_{\ell-k+1}) v_1\rangle e^{-\langle \mu, (I-\Ref_{\ell-k+1}) x\rangle}\right].
\end{eqnarray*}
Equating the coefficients with (\ref{eq:tildepid}), we find that $\pi^\mu$ satisfies (\ref{pde})--(\ref{bc2})
if for $k=1,\ldots,\ell$,
\[
c_{\ell-k} \langle \mu, (I-\Rot_{\ell-k})v_1\rangle\langle \mu, (I-\Ref_{\ell-k+1}) v_2\rangle=c_{\ell-k+1} \langle \mu, (I-\Ref_{\ell-k+1}) v_1\rangle\langle \mu, (I-\Rot_{\ell-k})v_2\rangle,
\]
and this is the recursion given in Lemma \ref{lem:vandermonde} (the condition $\t_\mu \notin \Theta_\ell$ guarantees that none of the four inner products in \eqref{eq:recck} is zero).

Now that we know that $\pi^\mu$ satisfies (\ref{pde})--(\ref{bc2}), it remains to show that it is a multiple of the stationary density.
Proposition~\ref{prop:nosignchange} and its proof show that $\pi^\mu$ is integrable and single signed.
Moreover, Proposition~\ref{prop:corner} implies that $\pi^\mu(0)=0$.
We therefore conclude from Proposition~\ref{prop:partdiffeq} that (BAR) is satisfied with $d\nu_0=\pi^\mu \, d\sigma$ and $d\nu_i=\pi^\mu/2 \, d\sigma_i$.
The claim thus follows from Proposition~\ref{prop:BAR}.

\section*{Acknowledgements}
We are grateful to Neil O'Connell for drawing our attention to \eqref{eq:duality}, which
motivated our research, and for several discussions. We would also like to thank Ivo Adan,
Gerry Foschini, Mike Harrison, and Ruth Williams for stimulating discussions. This work was
financially supported in part by the Science Foundation Ireland, grant number
SFI04/RP1/I512.

{\small
\bibliography{../../../bibdb}
\bibliographystyle{amsplain}
}
\end{document}